\documentclass{svproc}
\usepackage{amsmath}
\usepackage{xcolor}
\usepackage{amssymb}
\usepackage{enumitem}
\usepackage{float}
\usepackage{graphicx}
\usepackage{url}
\usepackage[square,numbers]{natbib}
\usepackage[capitalise]{cleveref}

\begin{document}
\mainmatter              
\title{A Novel Use of Pseudospectra in Mathematical Biology: Understanding HPA Axis Sensitivity}
\titlerunning{Pseudospectra for the HPA Axis}  
%
\author{Catherine Drysdale\inst{1} \and Matthew Colbrook\inst{2}}
\authorrunning{Drysdale and Colbrook} 
%
\tocauthor{Catherine Drysdale, Matthew Colbrook}
\institute{SMQB, University of Birmingham, Edgbaston, Birmingham, B15 2TT
\email{c.n.d.drysdale@bham.ac.uk},\\
\and
DAMTP, University of Cambridge, Wilberforce Rd, Cambridge, CB3 0WA}

\maketitle              

\begin{abstract}
The Hypothalamic-Pituitary-Adrenal (HPA) axis is a major neuroendocrine system, and its dysregulation is implicated in various diseases. This system also presents interesting mathematical challenges for modeling. We consider a nonlinear delay differential equation model and calculate pseudospectra of three different linearizations: a time-dependent Jacobian, linearization around the limit cycle, and dynamic mode decomposition (DMD) analysis of Koopman operators (global linearization). The time-dependent Jacobian provided insight into experimental phenomena, explaining why rats respond differently to perturbations during corticosterone secretion's upward versus downward slopes. We developed new mathematical techniques for the other two linearizations to calculate pseudospectra on Banach spaces and apply DMD to delay differential equations, respectively. These methods helped establish local and global limit cycle stability and study transients. Additionally, we discuss using pseudospectra to substantiate the model in experimental contexts and establish bio-variability via data-driven methods. This work is the first to utilize pseudospectra to explore the HPA axis.

\end{abstract}
\section{Introduction}

\begin{figure}[t]
\centering
\includegraphics[scale = 0.35]{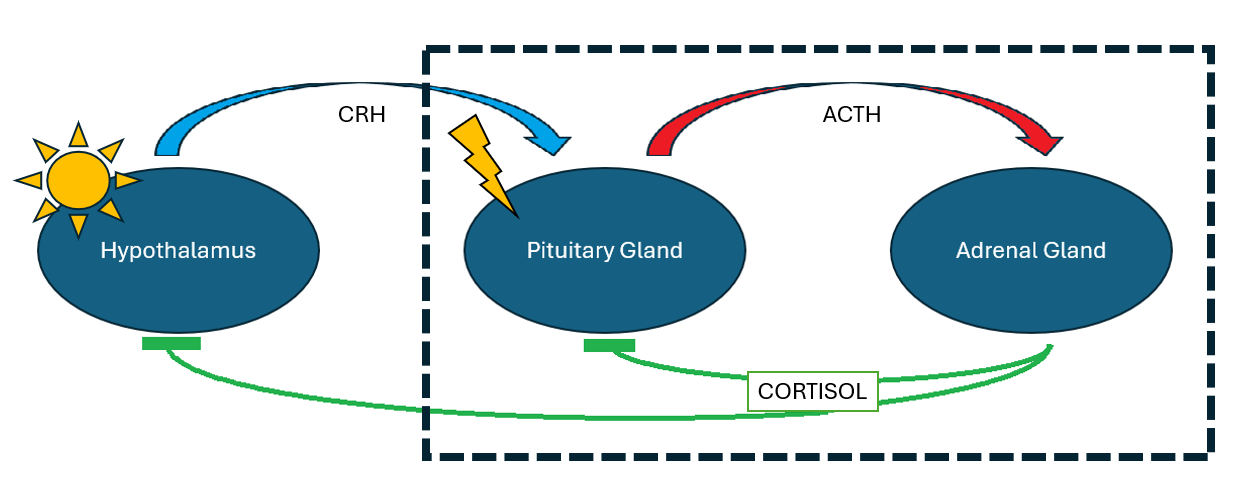}
\caption{
A schematic of the HPA axis. The sun represents that CRH is entrained by the day and night cycle, which triggers an enveloping of cortisol and ACTH on a circadian level. On an ultradian level, ACTH triggers cortisol, which then inhibits ACTH. Both operate on regular pulses, the amplitude of which varies according to the other and the wider system. The lightning bolt represents external stressors that cause ACTH and, thereby, cortisol to rise. The dashed line represents the restriction we make to the negative feedback loop between ACTH and cortisol that we consider in this paper.}
\label{fig1}
\end{figure}

\raggedbottom
The Hypothalamic-Pituitary-Adrenal (HPA) axis, a major neuroendocrine system, comprises the hypothalamus, pituitary gland, and adrenal glands. The hypothalamus releases corticotropin-releasing hormone (CRH), which stimulates the pituitary gland to secrete adrenocorticotropic hormone (ACTH). ACTH prompts the adrenal glands to produce cortisol, completing a negative feedback loop that inhibits further secretion of ACTH and CRH. These hormones follow circadian and ultradian rhythms. Three interacting behaviors of the HPA axis are noted: the circadian rhythm initiated by the hypothalamic suprachiasmatic nucleus (SCN) through CRH production, the ultradian pulsatility of ACTH and cortisol -- arguably the axis's intrinsic dynamics as they persist without CRH, and activation by external stressors, which increase hormone amplitudes but not frequency. \cref{fig1} illustrates this negative feedback loop. Dysregulation of the HPA axis, such as chronic over-activation resulting in sustained high cortisol levels or insufficient cortisol production, is implicated in various diseases, including anxiety disorders \cite{hpaanxietyreview}, depression \cite{hpadepressionreview}, PTSD \cite{hpaptsd}, Addison's disease \cite{hpaaddisons}, Cushing's syndrome \cite{hpacushing}, and chronic fatigue syndrome \cite{hpachronicfatigue}.

Mathematical modeling of the HPA axis has varied over the years, often targeting specific diseases or aligning with experimental observations. These models are typically nonlinear ordinary differential equations (ODEs) or delay differential equations (DDEs). Multiple mechanisms, such as time delays, negative feedback loops, nonlinearities, and varying timescales, can induce oscillations, as overviewed in \cite{novak2008design}. The ``minimal model" \cite{vinther2011minimal} incorporates CRH, ACTH, and cortisol and assumes a differentiable negative feedback effect on ACTH and CRH. Under minimal assumptions, this model demonstrates oscillations through a Hopf bifurcation, albeit with unphysical parameters. Another model \cite{depressedmodel} features similar components but attributes oscillations to Hill-function type nonlinearities, noting the absence of oscillations with physiologically reasonable parameters. Interestingly, perturbations in this model led to two additional fixed points, indicative of ``hypercortisolemic" and ``hypocortisolemic" depressive groups, suggesting that mechanisms such as delays, rather than nonlinearity, might drive oscillations. This hypothesis is supported by \cite{walker1}, where a model including ACTH, cortisol, and glucocorticoid receptors demonstrated self-sustained oscillations with sufficient delay, despite CRH being modeled as a constant input.

While a normal form adeptly captures the bifurcation structure, it often fails to align well with empirical data. This issue may persist across the models above, regardless of the chosen model's correctness. The tension arises from the interaction between non-normality and nonlinearity, where existing tools focus either on one aspect or the other but seldom both. For example, studies in \cite{schmidhenningson,chomaz} have shown that a non-normal linear operator, despite having stable eigenvalues, can give rise to a wave packet that interacts with the nonlinearity. Such interactions render first-order approximations, such as normal forms that depict bifurcation phenomena, not representative of the entire solution. Moreover, relying exclusively on linear tools, such as spectra and pseudospectra, can fail to capture the transient behaviors of the complete nonlinear system.

We seek a balanced approach by exploring different linearizations that help us comprehensively understand the complete nonlinear system. We explore transient behaviors for nonlinear systems by computing different types of pseudospectra for various linearizations, gaining different insights from each linearization.

We consider the following  nonlinear DDE model of the HPA axis:
\begin{equation}\label{main_HPA}
    \frac{dA}{dt} = - e_a A +  \frac{hc^{m_1}}{c^{m_1} + [C(t-\tau_1)]^{m_1}},\quad
    \frac{dC}{dt} = - e_c C +  \frac{\beta [A(t-\tau_2)]^{m_2}}{a^{m_2} + [A(t-\tau_2)]^{m_2}},
\end{equation}
where $A$ denotes ACTH and $C$ denotes cortisol. The initial conditions and constants are described in \cref{sec_model}. This model, used in \cite{malek1, malek2}, where analysis demonstrates the existence of a Hopf bifurcation, was chosen for its ability to exhibit oscillatory phenomena at physiologically reasonable parameter values. It enables focused analysis of the intrinsic dynamics of the HPA axis, particularly the ultradian rhythms of ACTH and cortisol. We seek to answer the question, ``What should we linearize around?" and explore three linearizations:
\begin{itemize}[label=\textbullet]
    \item A time-dependent Jacobian with time-dependent spectra;
    \item Linearization around the limit cycle, investigating whether there is inherent non-normality that could suggest transient effects;
    \item A dynamic mode decomposition (DMD) approach, which transforms a finite-dimensional nonlinear system (that approximates our original DDE system) into an infinite-dimensional linear system (with no delays). We apply the data-driven method ResDMD \cite{colbrook2021rigorousKoop} to compute pseudospectra.
\end{itemize}

To our knowledge, this work is the first to use pseudospectra to explore the HPA axis. Different linearizations allow us to investigate different biological phenomena. For instance, by considering the time-dependent Jacobian, we see that there is instability corresponding to the upward slopes of cortisol, as in rat models when they are seen to become more aggressive if presented with a stressor on the upward slope of corticosterone\footnote{Corticosterone is the equivalent of cortisol in rats. We use mathematical models fitted to human data. Hence, cortisol best describes the variable $C$ in \cref{main_HPA}. However, relevant experiments were done on rats, so we often focus on the similarities between the behavior of cortisol in our models and corticosterone in rats.} \cite{rats1, rats2}. On the downward slope, the non-normality index is low, suggesting the system is not sensitive to perturbation during this phase. This is the first mathematical modeling to match the experimental findings. 

The paper is organized as follows. \cref{sec_model} provides a detailed description of the model, including the motivations for the parameters and the process of de-dimensionalizing the model. \cref{sec_pseudospec}, titled ``pseudospectra of the model three ways", presents the pseudospectra calculated from a time-dependent Jacobian, linearization around the limit cycle, and DMD. In the conclusion, we discuss our results and the contexts in which each can be useful for explaining various experiments and observations. Additionally, we explore how pseudospectra might be used for model validation, i.e., determining whether the perturbed model captures experimental results. We also consider the potential of using data-driven methods to capture the spectrum and pseudospectrum. Code producing all of the figures of the paper can be found here: \textcolor[rgb]{0,0,1}{\url{https://github.com/miinadietrich/hpa_axis_pseudospectra}}.

\section{Model}\label{sec_model}
\begin{figure}[t]
\centering
\includegraphics[height=4.2cm,trim={0mm 0mm 0mm 0mm},clip]{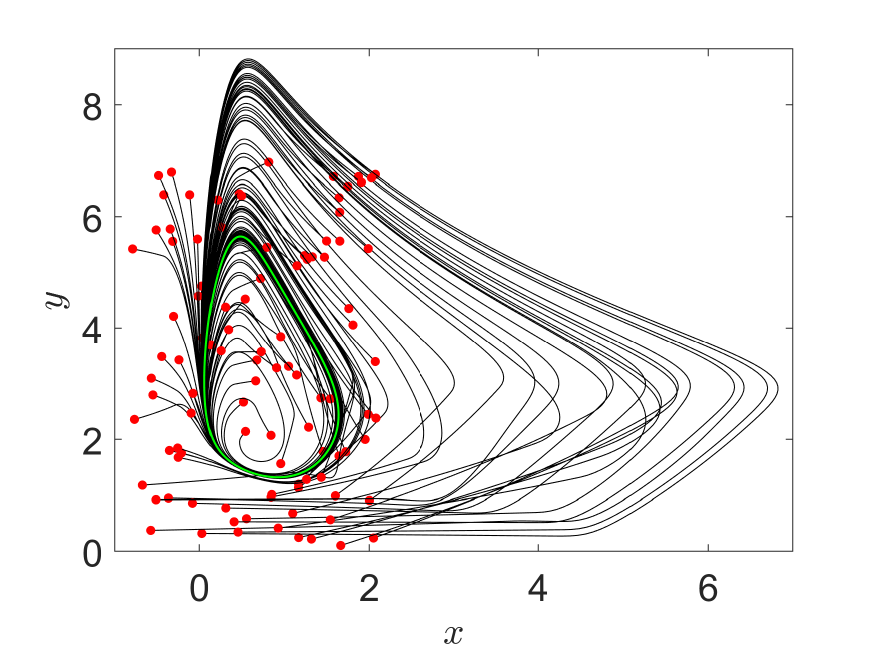}
\hspace{5mm}
\includegraphics[height=3.9cm,trim={0mm 0mm 0mm 0mm},clip]{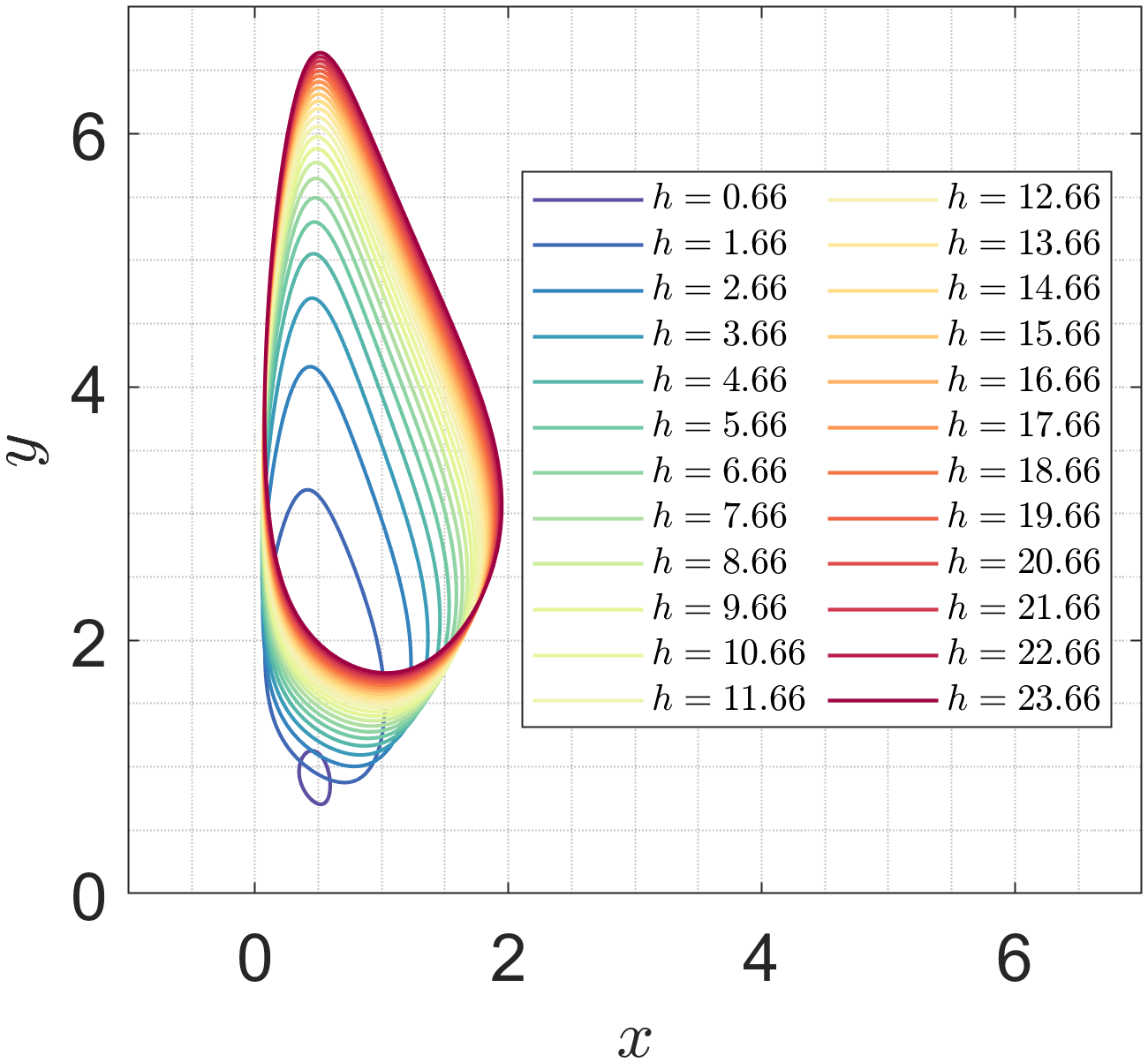}\vspace{-4mm}
\caption{Left: Trajectories with 100 randomly selected constant initial conditions (red dots) for $h=7.66$ and limit cycle (green). Right: The limit cycle as $h$ is varied.}
\label{fig_curve}
\end{figure}

We use the following constants from \cite{malek1, malek2} for the system in \cref{main_HPA}:
\begin{itemize}
    \item $e_a = 0.04$, the rate of elimination of ACTH;
    \item $e_c = 0.01$, the rate of elimination of cortisol;
    \item $m_1 = 4$, the number of cortisol modules that bind to the free receptors in the adrenal gland;
    \item $m_2 = 4$, the number of ACTH modules that bind to the receptors in the pituitary gland;
    \item $a = 21$, the Hill half-maximum constant for ACTH;
    \item $c= 6.11$, the Hill half-maximum constant for cortisol; 
    \item $h =7.66$, the positive stimulation of CRH on ACTH secretion;
    \item $\beta = 1$, the rate of production of cortisol from the release of ACTH;
    \item $\tau_1 = 15$, the delay for secreted cortisol to have negative feedback on ACTH;
    \item $\tau_2 = 15$, the delay for cortisol to stimulate ACTH production in the adrenal glands.
\end{itemize} We non-dimensionalize by letting $A = ax$ and $C = cy$ and $t = \tau/e_c$. We also scale the delays accordingly with $t_1 = e_c\tau_1$ and $t_2 = e_c\tau_2$ and set $c_1=e_a/e_c=4$, $c_2=1/(ae_c)\approx 4.76$, $c_3=\beta/(c e_c)\approx 16.37.$ The result is our system of study:
\begin{equation}\label{our_equation}
    \frac{dx}{d\tau} = - c_1 x +  \frac{h  c_2}{1+ [y(\tau-t_1)]^{m_1}} ,\quad
    \frac{dy}{d\tau} = - y + c_3 \frac{[x(\tau - t_2)]^{m_2}}{1 + [x(\tau-t_2)]^{m_2}}
\end{equation} 
with initial conditions $x(t) = x_0$ for $ t \in [-t_1, 0]$ and $y(t) = y_0$ for $ t \in [-t_2, 0]$.
A default initial condition is $(x_0,y_0) = ( 0.8858, 1.7461)$, arising from solving \cref{main_HPA} for a fixed point by linearizing and setting the left-hand-side equal to zero \cite{malek1}. In the same vein as the bifurcation analysis, we also consider a range of physical $h$ values, i.e., for different CRH inputs that would occur throughout the day as CRH operates on a circadian cycle. In \cite{malek2}, the default value of $h$ was $7.66$, but the full range was from $3.068$ to $23$. For the considered range of values of $h$, the full system has a limit cycle behavior that we plot in \cref{fig_curve}. This limit cycle varies continuously with $h$.

\section{Pseudospectra of the model three ways}\label{sec_pseudospec}

We now consider three different families of pseudospectra associated with the above model. Readers from mathematical biology may be more familiar with spectra or eigenvalues, which are used to study a system's long-term stability (or instability). In other words, spectra capture the \textit{post-transient} behavior of linear systems. In contrast, pseudospectra enable us to quantify sensitivity to perturbations and capture the \textit{transient} dynamics of linear systems. An excellent account can be found in \cite[Chapters 14-19]{trefethen}.

For instance, consider a finite matrix $A$ whose eigenvalues all have magnitude less than $1$. The linear discrete-time dynamical system $x_{k+1} = Ax_k$ is stable since $\lim_{k\rightarrow\infty}\|A^{k}\|\rightarrow0$. However, the pseudospectral contours may protrude beyond the unit circle, indicating that there could be some transient growth before decay when the system is perturbed. In the case of pseudospectra calculated on a limit cycle, this could imply that the limit cycle takes longer to stabilize after a perturbation. We will discuss this point below and compute so-called Kreiss constants that capture this behavior.

Before looking at the different types of pseudospectra of linearizations of \cref{our_equation}, it is instructive to consider the following three equivalent definitions of pseudospectra of finite matrices.

\begin{definition}[Pseudospectra of Matrices]\label{pseudospectra_def}
 Let $A\in \mathbb{C}^{N \times N}$ and $\epsilon > 0 $. The $\epsilon$-pseudospectrum $\mathrm{Sp}_{\epsilon}(A)$ is the set of $z \in \mathbb{C}$ defined equivalently by any of: 
    \begin{itemize}
        \item[$\bullet$] 
        $\|(A-z)^{-1}\|^{-1} \leq \epsilon $
        \item[$\bullet$] 
        $ z \in \mathrm{Sp}(A+E)$ for some matrix $E\in\mathbb{C}^{N\times N}$ with $\|E\|\leq\epsilon$
        \item[$\bullet$] 
        $\|(A-z)u\| \leq \epsilon$ for some vector $u\in\mathbb{C}^{N}$ with $\|u\| = 1$.
    \end{itemize}
\end{definition}
We can see from the first point in the definition that eigenvalues are contained within the $\epsilon$-contours of the function $z\mapsto\|(A-z)^{-1}\|^{-1}$. It is often an excellent numerical check to ensure that spectral computations are accurate by verifying that the spectrum is contained within the pseudospectrum. The second condition describes how far the eigenvalues move in response to a size $\epsilon$ perturbation. The third condition associates spatial structures (i.e., the vector $u$) with the perturbed eigenvalues, which often have physical meaning.

Often, for nonlinear systems, there is a choice regarding how to obtain the linearization or the matrix $A$ in the above definitions. This is more nuanced in the context of delay-differential equations. Furthermore, the $\epsilon$-pseudospectra form distinct sets; therefore, comparing the pseudospectra of one operator to another is also nuanced. In each section, we provide a short preamble explaining what we hope to understand from this linearization and the method used with numerical details. Then, we present the results and their biological ramifications.

\subsection{Time-dependent Jacobian: Pointwise linearization}

We linearize \cref{our_equation} with a time-dependent Jacobian around an arbitrary time point rather than an equilibrium point. This approach is often used in feedback control \cite[Chapter 12]{khalil}. It can be viewed as representing a nonlinear system by a linear system over an instantaneous time window.

\subsubsection{Method:}
Taking a base trajectory $(x_0,y_0)$, we write
$
x(\tau) = x_0(\tau) + \tilde{x}(\tau)
$
and
$
y(\tau)= y_0(\tau) + \tilde{y}(\tau).
$
To first order, this yields the following linear DDE system:
\begin{equation}
\label{pseudospec_eq1}
      \begin{pmatrix}
    \frac{d\tilde{x}}{d\tau} \\ \frac{d\tilde{y}}{d\tau} \end{pmatrix} =  A
    \begin{pmatrix} \tilde{x} \\ \tilde{y} \end{pmatrix}
    +   B \begin{pmatrix} \tilde{x}(\tau - t_1) \\ \tilde{y}(\tau - t_1) \end{pmatrix}
    +   C \begin{pmatrix} \tilde{x}(\tau - t_2) \\ \tilde{y}(\tau - t_2) \end{pmatrix},
\end{equation}
where $A,B$ and $C$ are the following matrices:
$$
    A = \begin{pmatrix} -c_1 & 0 \\ 0 & -1 \end{pmatrix},\quad\!\!\!
    B = \begin{pmatrix} 0 & - hc_2 \frac{m_1 [y_0(\tau - t_1)]^{m_1 -1}}{(1 +[y_0(\tau - t_1)]^{m_1})^2}\\
    0 & 0\end{pmatrix},\quad\!\!\!
    C= \begin{pmatrix} 0 & 0 \\ c_3 \frac{m_2 [x_0(\tau - t_2)]^{m_2 -1}}{(1 + [x_0(\tau - t_2)]^{m_2})^2} & 0 \end{pmatrix}. 
$$
Taking Laplace transforms of \eqref{pseudospec_eq1} leads to 
$
    \Delta_\tau(\lambda) = \lambda I - A - Be^{-\lambda t_1} - C e^{-\lambda t_2}.
$
In particular, we obtain a time-dependent spectrum
$$
\mathrm{Sp}(\Delta_\tau)=\{z\in\mathbb{C}:\mathrm{det}(\Delta_\tau( z)) = 0\}
$$
corresponding to a \textit{nonlinear} eigenvalue problem \cite{guttel2017nonlinear}. The spectrum is generally infinite (in contrast to linear finite matrix eigenvalue problems), but for any $c\in\mathbb{R}$, there are only finitely many eigenvalues $\lambda$ with $\mathrm{Re}(\lambda)\geq c$ \cite[Theorem 1.5]{gu2003stability}. In what follows, we use the MATLAB package TDS-CONTROL \cite{appeltans2022tds} to compute the (finitely many) eigenvalues in a right-half plane.

We are interested in pseudospectra, defined as
$$
\mathrm{Sp}_\epsilon(\Delta_\tau)=\{z\in\mathbb{C}:\exists u\in\mathbb{C}^2\text{ s.t. }\|\Delta_\tau(z)u\|\leq\epsilon,\|u\|=1\}.
$$
Note how this generalizes the matrix case in \cref{pseudospectra_def}. One can show that the other two equivalent definitions in \cref{pseudospectra_def} also hold.\footnote{For the second point in \cref{pseudospectra_def}, we consider pencils $\lambda\rightarrow E(\lambda)$ bounded by $\epsilon$.} To compute $\mathrm{Sp}_\epsilon(\Delta_\tau)$, we compute the smallest singular values of $\Delta_\tau(z)$ over a grid of $z$ points.

We are interested in how the eigenvalues and pseudospectra vary with time $\tau$ (for the selection of base points $(x_0, y_0)$) and the parameter $h$ to see if they match the experimental findings of more transient growth regarding upward peaks versus downward peaks of cortisol. We consider
the spectral abscissa, defined as
$
\alpha=\alpha(\tau) = \max_{\lambda\in\mathrm{Sp}(\Delta_\tau)}\mathrm{Re}(\lambda),
$
and the distance to instability, defined as
$
d_\tau=\inf\{\epsilon>0:\exists s\in\mathbb{R},\|\Delta_\tau( is)^{-1}\|^{-1}=\epsilon\}
$
whenever $\alpha(\tau)<0$. The system in \cref{pseudospec_eq1} is (asymptotically) stable if and only if $\alpha(\tau)<0$ \cite[Chapter 1]{gu2003stability}. In that case, $d_\tau$ tells us the minimum size of perturbation needed to cause the system to become unstable. We also consider the ratio $-\alpha/d_\tau$, whenever $\alpha<0$. This index tells us the ratio between the distance of the spectrum to the imaginary axis and the minimum value of $\|\Delta_\tau( z)^{-1}\|^{-1}$ on the imaginary axis. With an abuse of terminology, we shall refer to this as a non-normality index.\footnote{Since $\Delta_\tau(\lambda)$ is nonlinear in the spectral parameter $\lambda$, the phrase non-normality here is used in the sense of large pseudospectral regions enclosing eigenvalues.}

\subsubsection{Results:} 
\begin{figure}[t]
\centering
\includegraphics[width=0.8\textwidth,trim={0mm 0mm 0mm 0mm},clip]{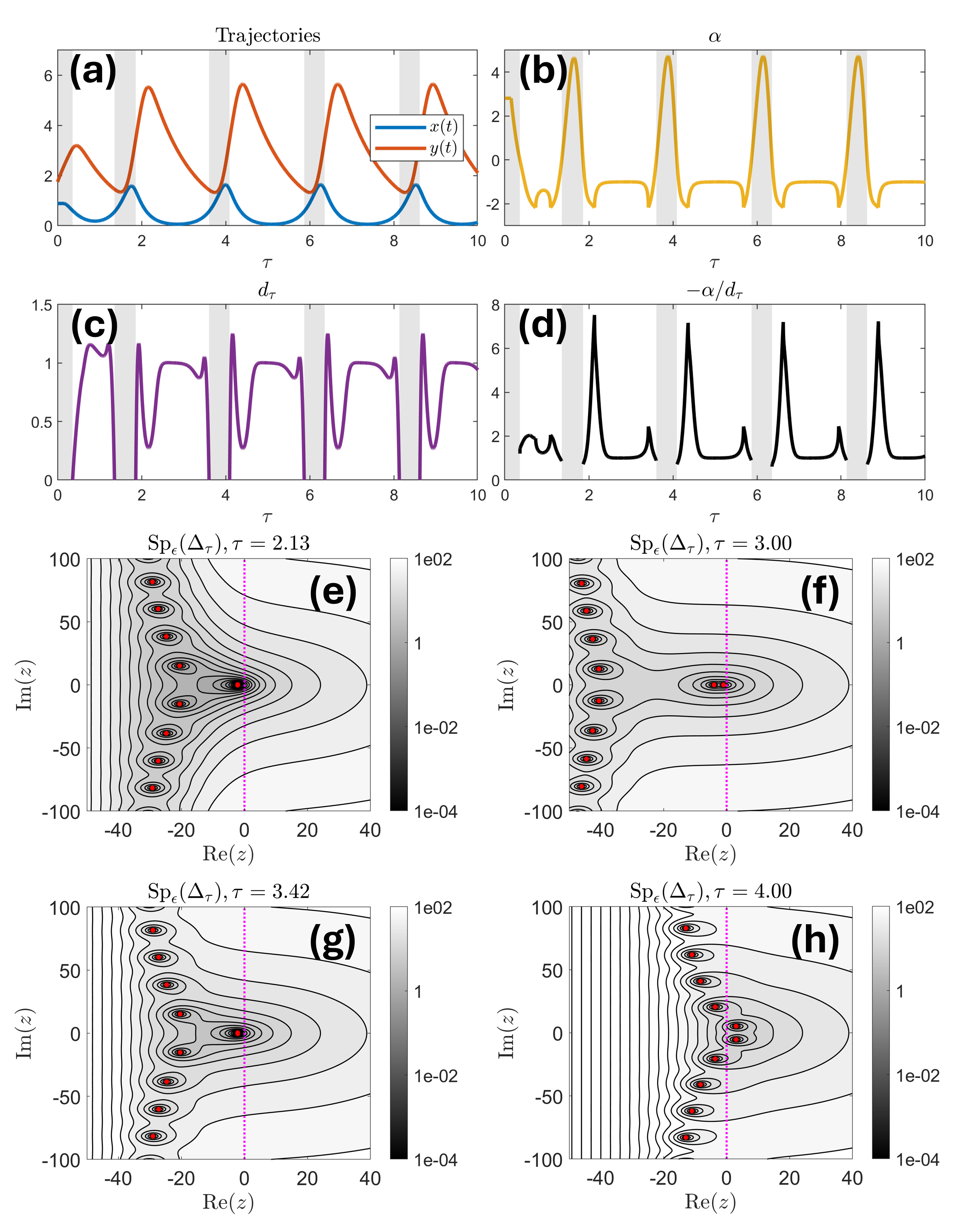}\vspace{-5mm}
\caption{(a) -- (d): Solution trajectories, spectral abscissa, distance to instability, and ratio $-\alpha/d_\tau$. (e) -- (h): Pseudospectra of the linearized system for different $\tau$. We display pseudospectra by plotting contours of $\epsilon$ (colorbars). The eigenvalues are shown in red. The stability line (imaginary axis) is the dotted magenta line.}\vspace{-5mm}
\label{default_h_fig}
\end{figure}

\cref{default_h_fig} (panels (a) -- (d)) plots $\alpha$, $d_\tau$, and $-\alpha/d_\tau$ for $h=7.66$. We also show pseudospectra plots at different time points (panels (e) -- (h)).

The instability regions correspond to the upward rise of cortisol and ACTH (see panel (a)). They also coincide with points where $x$ (ACTH) is maximized, and $y$ (cortisol) is minimized before it increases again. These observations strongly correlate with the findings in \cite{rats1, rats2}, where rats became more aggressive under a stressor during the upward slope of corticosterone rather than during the downward peak. Additionally, the downward slope of cortisol is entirely within the stable region of the graph.

The index $-\alpha_\tau/d_\tau$ has peaks before and after the unstable region (panel (d)). These correspond to the dips in the distance to instability (panel (c)) and the spectral abscissa (panel (b)). Therefore, a small perturbation can easily bring us into the unstable half-plane, and this in itself should be considered as a ``non-normal feature". In \cite{lightman2020dynamics}, it was said that the HPA axis needs to be both ``sensitive to environmental perturbations, and able to respond differently to both small and large stimuli", as well as be ``robust with preservation of dynamic behavior during these perturbations". We can associate the two different-sized peaks with different responses to small and large stimuli at different times. The fact that there is a plateau between these two peaks possibly corresponds to robustness and the system's desire to return to equilibrium following a perturbation but being sufficiently ``non-normal" at other points to respond to perturbation. These ideas have to be experimentally validated.

The leading eigenvalue splits into two eigenvalues (panel (e) to panel (f)) on the stable side of the plane. These eigenvalues coalesce again (panel (g)) before splitting and going into the right-half plane (panel (h)). The coalescence corresponds to the two unstable peaks we saw before, whereas the separation corresponds to the instability. Arguably, we have an ``intermittent" Hopf bifurcation. Such bifurcations have been used for the relation between CRH drive and delays \cite{zavala2019mathematical}, but not so far in the study of the ultradian rhythms.  

\begin{figure}[t]
    \centering
    \includegraphics[width=0.4\textwidth,trim={0mm 0mm 0mm 0mm},clip]{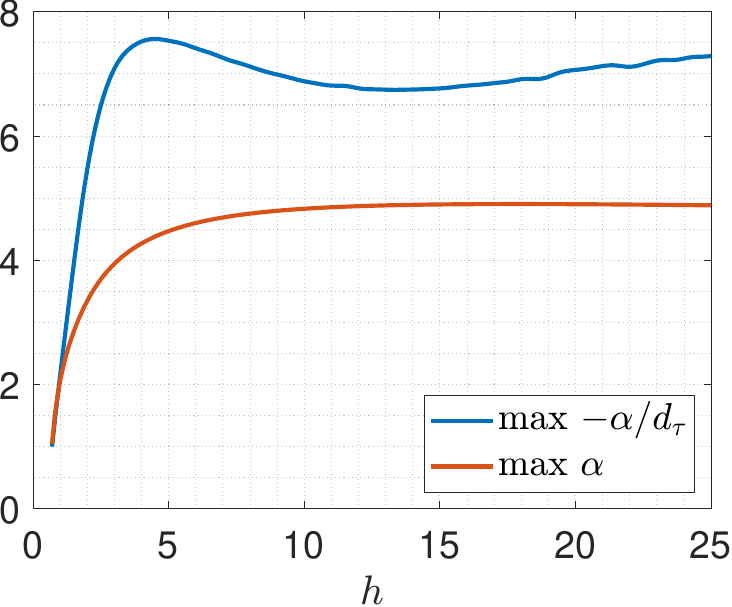}
    \caption{Maximum $\alpha$ and index $-\alpha/d\tau$ (over the limit cycle) as $h$ varies.}
    \label{vary_h_fig}
    \end{figure}

This behavior is consistent across different values of $h$, i.e., independent of CRH input. \cref{vary_h_fig} plots the maximum index $-\alpha/d_\tau$ and $\alpha$ as $h$ varies.

\subsection{Floquet pseudospectrum: Linearization about limit cycle}
\label{sec:floquet}

We now consider a periodic linear DDE system, corresponding to linearization about the limit cycle. Unlike linear periodic ODE systems, a complete Floquet theory does not generally exist. Nevertheless, we can still characterize stability and transients through the spectral properties of an infinite-dimensional operator.\footnote{In our case, this operator is completely continuous. This implies that its spectrum is countable, zero is the only accumulation point of the spectrum, and all non-zero parts of the spectrum are eigenvalues of finite multiplicity.}

\subsubsection{Method:} If the base trajectory $(x_0,y_0)$ is periodic, then \cref{pseudospec_eq1} is a linear periodic DDE system, where the matrices $A$, $B$ and $C$ are periodic in time. The stability of such systems is characterized by characteristic multipliers \cite[Chapter 8]{Hale}, which are similar to those found in Floquet theory.\footnote{It should be noted that a complete Floquet theory does not generally exist.} Let $\omega$ denote the period of the limit cycle. Consider the operator $U$, which maps a continuous function $\phi$ defined over the interval $[-t_1,0]$ to the solution of \cref{pseudospec_eq1} over the interval $[-t_1+\omega, \omega]$, starting from the initial condition $\phi$. The operator $U$ acts on the Banach space $C([-t_1,0])$ of continuous functions on $[-t_1,0]$ equipped with the $L^\infty$ norm. The stability and transient behavior of the system are determined by the spectra and pseudospectra of $U$. Pseudospectra of $U$ are given by\footnote{A subtlety arises here. The spectrum of $U$ is independent of the initial time for which we study the system in \cref{pseudospec_eq1}, and the eigenspaces for different initial times are diffeomorphic. However, the pseudospectra can change. For consistency, as we vary $h$, we take a point at which ACTH peaks as the initial time.}
$$
\mathrm{Sp}_\epsilon(U)=\{z\in\mathbb{C}:\exists\phi\in C([-t_1,0])\text{ s.t. }\|(U-zI)\phi\|_{L^\infty}\leq\epsilon,\|\phi\|_{L^\infty}=1\}.
$$
Computing pseudospectra in Banach spaces generally poses a significant challenge but is feasible in our case. Consider a discrete grid of equally spaced points $s_1, s_2, \ldots, s_N$ with $s_1 = -t_1 < s_2 < \cdots < s_N = 0$. We utilize piecewise affine hat functions, $\phi_j$, such that $\phi_j(s_i) = \delta_{ij}$. These functions form a basis of the piecewise affine functions on the interval $[-t_1,0]$ with knots at the points $s_j$. Using this basis for the variables $x$ and $y$, we compute $U[\phi_{j_1}; \phi_{j_2}]$ by solving the periodic DDE system and then evaluate the solution at the same grid points to form $T\in\mathbb{C}^{2N\times 2N}$ which approximates $U$. For a given $z$, we compute the inverse $(T-zI)^{-1}$ and subsequently calculate its matrix $l^\infty$ norm. This provides an approximation of $\|(U-zI)^{-1}\|^{-1}$, which characterizes the level sets of $\mathrm{Sp}_\epsilon(U)$.

We also consider the Kreiss constant of $U$, 
$
\mathcal{K}_c(U)=\sup_{\epsilon>0}[\sup\{|z|:z\in\mathrm{Sp}_\epsilon(U)\}-c]/\epsilon=\sup_{|z|>c}(|z|-c)\|(U-z)^{-1}\|
$
defined for $c>\sup\{|z|:z\in\mathrm{Sp}_\epsilon(U)\}$. These constants give a lower bound on transient behavior about the limit cycle through the relation \cite[Chapter 16]{trefethen}:
$
\sup_{k\geq 0}c^{-k}\|U^k\|\geq \mathcal{K}_c(U).
$

\subsubsection{Results:} 
\begin{figure}[t]
    \centering
    \includegraphics[width=0.4\textwidth,trim={0mm 0mm 0mm 0mm},clip]{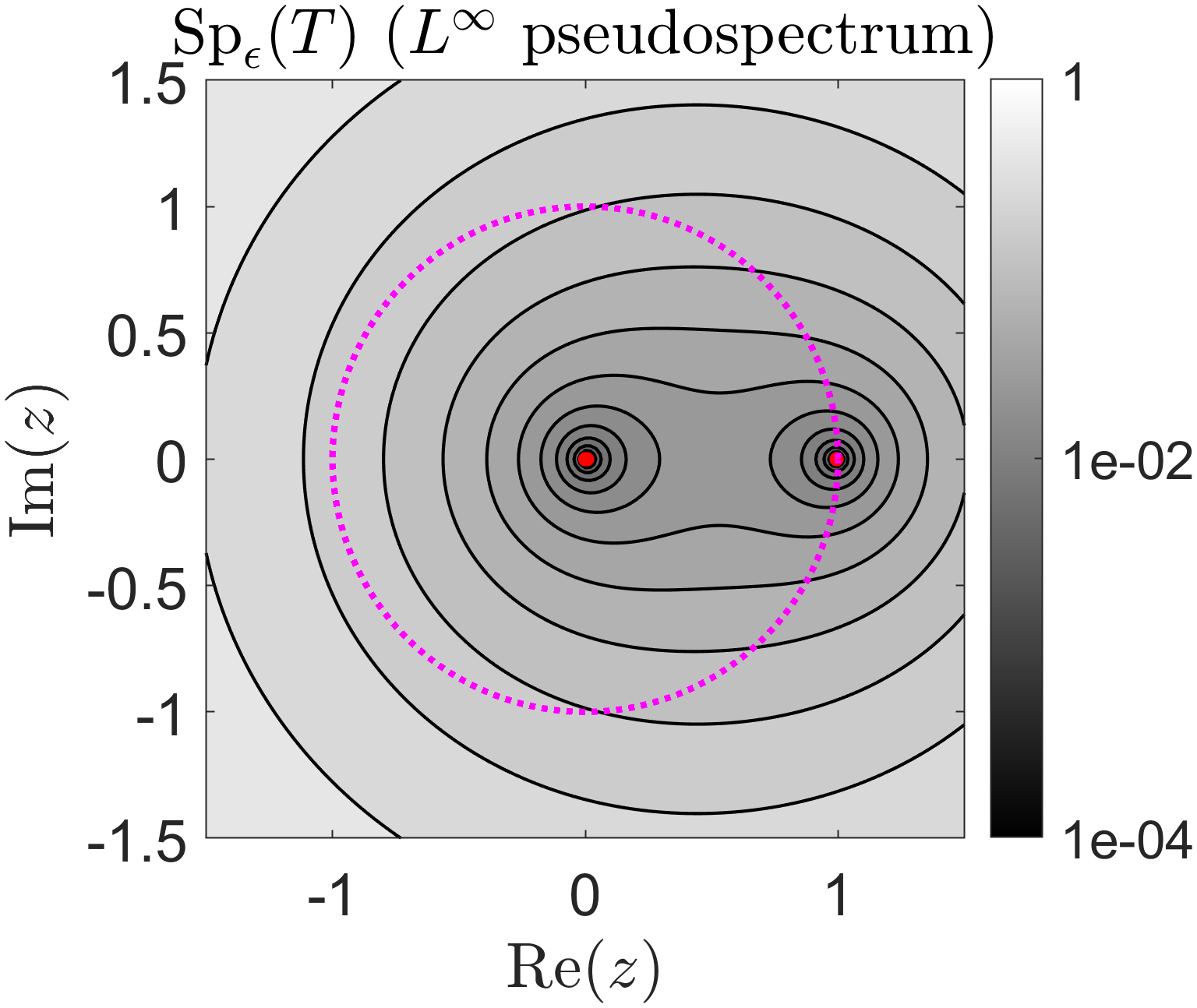}
    \hspace{10mm}\includegraphics[width=0.4\textwidth,trim={0mm 0mm 0mm 0mm},clip]{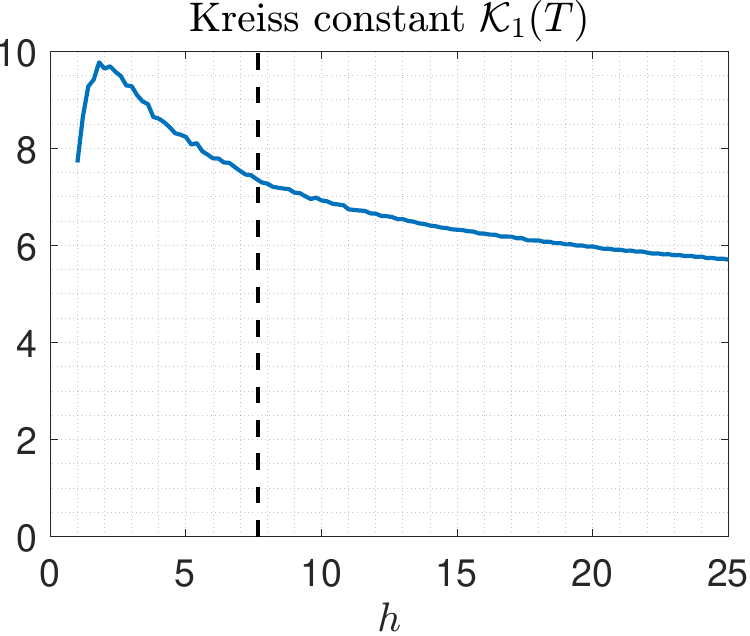}
    \caption{Left: The pseudospectrum of $T$ for the default parameter $h=7.66$. The eigenvalues are shown in red, and the stability circle (unit circle) is the dotted magenta circle. Right: The Kreiss constant for different values of $h$ (the dashed line is $h=7.66$).}
\label{vary_h_kreiss_fig}
\end{figure}

We take $N=50$, which corresponds to a total of $100$ basis functions. \cref{vary_h_kreiss_fig} (left) displays the pseudospectra for the default $h=7.66$. There is a dominant eigenvalue $\lambda\approx 0.9946$ (which has converged in $N$) very close to one, while all other eigenvalues cluster at zero; however, large regions of pseudospectra are present. We consider $c=1$ and found $\mathcal{K}_1(T)=7.4014$. \cref{vary_h_kreiss_fig} (right) plots the Kreiss constant for different values of $h$. After an initial peak, we see a decrease in transient effects with increasing $h$.

The limit cycle is stable, as expected in a hormone system. However, it may exhibit significant transient growth under perturbation, as exhibited by the large regions of pseudospectra. In particular, following a perturbation, the dynamics may take longer to stabilize compared to a system described by a normal matrix (whose pseudospectra are well-behaved). Given that the dominant eigenvalue is very close to one, these transient behaviors will likely be critically important. The variation in the Kreiss constant in \cref{vary_h_kreiss_fig} suggests that the limit cycle may be more susceptible to perturbations when CRH concentrations vary at different times of the day.

\subsection{DMD analysis: Global linearization}

In our final example, we consider a global linearization through Koopman operators. These operators act on a lifted infinite-dimensional space of functions of the state $(x,y)$. We use data-driven algorithms to compute a Galerkin approximation of this operator and its pseudospectra from trajectory data.

\subsubsection{Method:} We first augment the state $(x,y)$ to form the infinite state vector $\mathbf{x}(\tau)=(x(\tau),y(\tau),x(\tau-t_1),y(\tau-t_1),x(\tau-2t_1),y(\tau-2t_1),\ldots)$. Formally, this leads to an infinite system of nonlinear ODEs for $\mathbf{x}$. As an approximation, we truncate $\mathbf{x}$ to $d-1$ time lags so that $\mathbf{x}\in\mathbb{R}^{2d}$ ($d=10$ in what follows), leading to an approximate system
$
\frac{d\mathbf{x}}{d\tau} = \mathbf{f}(\mathbf{x}),
$
for some nonlinear function $\mathbf{f}:\mathbb{R}^{2d}\rightarrow \mathbb{R}^{2d}$.
 Letting $\omega$ denote the time period of the limit cycle (which depends on $h$), we sample this system with $\Delta \tau=\omega/10$, yielding a discrete-time dynamical system
\begin{equation}
\label{Koopman_system}
\mathbf{x}_{n+1}=\mathbf{F}(\mathbf{x}_n),\quad n=0,1,2,\ldots. 
\end{equation}
Given a suitable function (called an observable) $g:\mathbb{R}^{2d}\rightarrow \mathbb{C}$, the Koopman operator $K$ acts on $g$ by composition with $\mathbf{F}$: 
$$
[Kg](\mathbf{x}) = g(\mathbf{F}(\mathbf{x})).
$$
The key point is that $K$ is linear. We have traded the finite-dimensional nonlinear system in \cref{Koopman_system} for a linear operator $K$ that acts on an infinite-dimensional function space of observables. In general, the choice of function space matters, see \cite{mezic2020spectrum,colbrook2023multiverse}, but it is typical to take an $L^2$ space with respect to some measure. Since $g(\mathbf{x}_n)=[K^ng](\mathbf{x}_0)$, we see that $K$ characterizes the forward-time dynamics of the system and its pseudospectra capture transient behavior. In particular, observables $g$ associated with $\mathrm{Sp}_\epsilon(K)$ so that $\|(K-zI)g\|\leq \epsilon\|g\|$ are physically significant and correspond to approximate coherence since
$
\|(K^n-z^nI)g\|=\|g(\mathbf{x}_n)-g(\mathbf{x}_0)\|=\mathcal{O}(n\epsilon\|g\|).
$

Koopman operators have received considerable attention over the last decade; see \cite{brunton2021modern,colbrook2023multiverse}. As well as dealing with nonlinearity, the considerable advantage of this approach is that we do not have to know the map $\mathbf{F}$. One can also use experimental data. In what follows, we will build approximations of $K$ using trajectory data. We will employ the Residual DMD (ResDMD) algorithm, which computes spectra and pseudospectra of $K$ with rigorous convergence guarantees. A complete discussion of this algorithm is beyond the scope or space of this paper, but we refer the reader to \cite{colbrook2021rigorousKoop,colbrook2024limits}.

We collect data from $10^4$ initial points selected uniformly at random in $[-3,5]\times[-1,8]$. We run the system in \cref{our_equation} forward by $(d+1)\times\Delta\tau$ to generate an initial $\mathbf{x}_0$ as well as $\mathbf{x}_1$, and $\mathbf{x}_2$. In the language of DMD, the data consists of $M=2\times 10^4$ snapshots formed by length-three trajectories for each initial condition. As our basis of observable, we use $k$-means on the data to select $N=400$ centers in $\mathbb{R}^{2d}$, which we use with Gaussian radial basis functions.

\begin{figure}[t]
\centering
\includegraphics[width=0.4\textwidth,trim={0mm 0mm 0mm 0mm},clip]{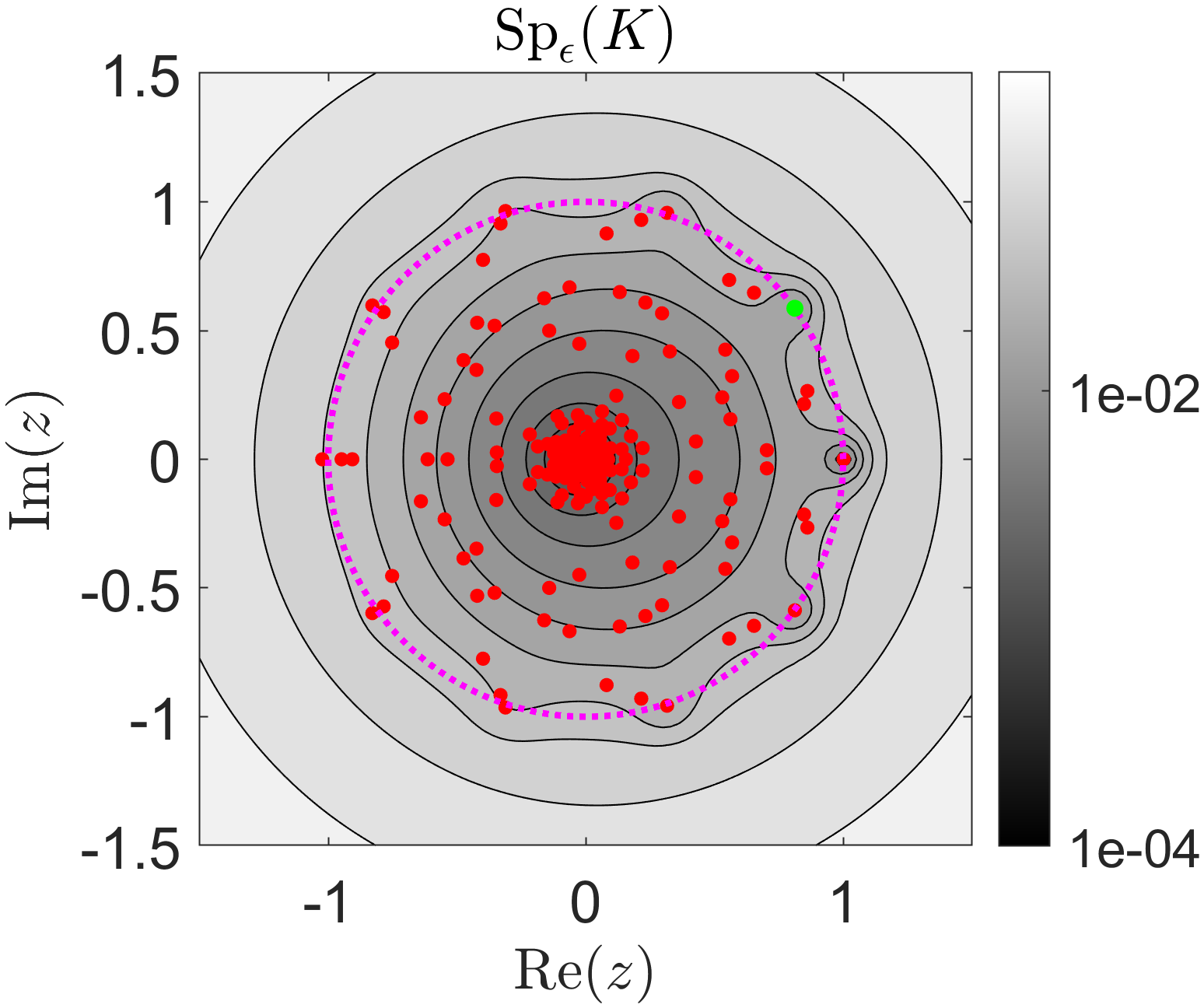}
\hspace{5mm}\includegraphics[width=0.4\textwidth,trim={0mm 0mm 0mm 0mm},clip]{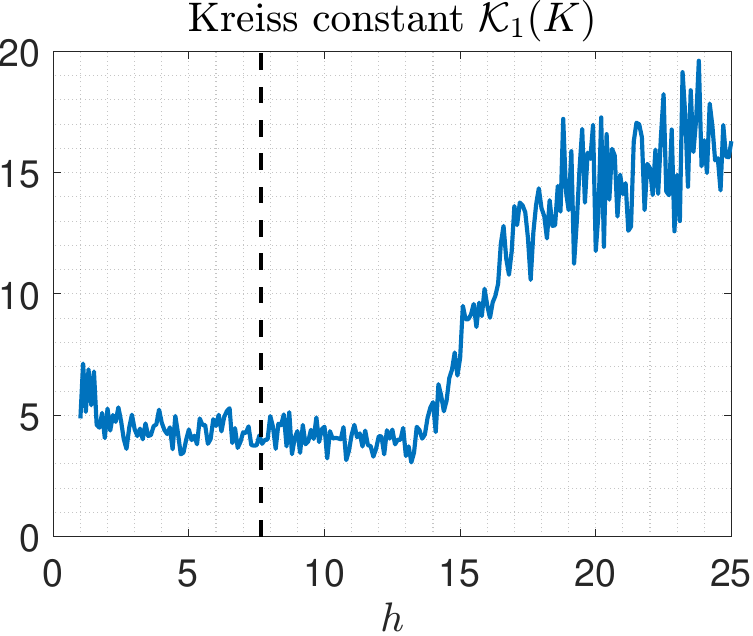}
    \caption{Left: The pseudospectrum of $K$ for the default parameter $h=7.66$. The eigenvalues are shown in red, and the stability circle (unit circle) is the dotted magenta circle. Right: The Kreiss constant for different values of $h$ (the dashed line is $h=7.66$).}
\label{DMD_kreiss_fig}
\end{figure}

\subsubsection{Results:} 
\begin{figure}[t]
    \centering
    \includegraphics[width=1\textwidth,trim={0mm 0mm 0mm 0mm},clip]{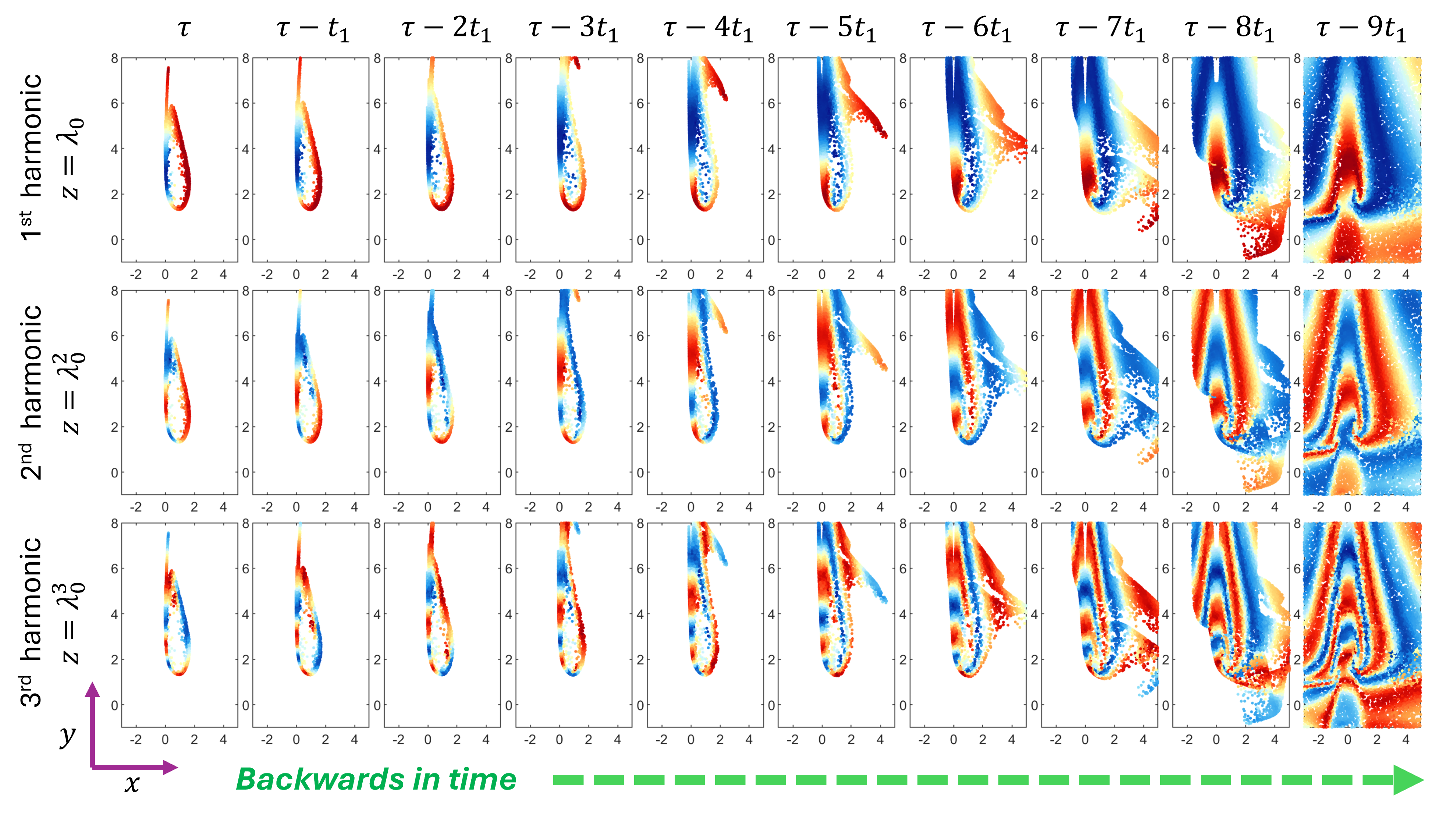}
    \caption{Eigenfunctions of $K$ corresponding to harmonics on the limit cycle (left). Moving from left to right, we march backward in time, capturing the transient dynamics off the limit cycle.}
\label{DMD_eigenfunction_fig}
\end{figure}

\cref{DMD_kreiss_fig} (left) shows the pseudospectrum computed by ResDMD, as well as the DMD eigenvalues (red dots). On the circle, a lattice of eigenvalues corresponds to harmonics on the limit cycle generated by an eigenfunction $g_0$ with eigenvalue $\lambda_0$ (green dot). The lattice structure arises since the Koopman operator is multiplicative: $K[g_1g_2]=[Kg_1][Kg_2]$. These harmonics are shown in \cref{DMD_eigenfunction_fig}, where we plot the eigenfunctions corresponding to $\lambda_0$, $\lambda_0^2$, and $\lambda_0^3$ across the entire state space $\mathbb{R}^{20}$. Here, the interpretation is that each two-dimensional block of $\mathbf{x}\in\mathbb{R}^{20}$ corresponds to marching the eigenfunction backward in time, capturing dynamics off the limit cycle.

\cref{DMD_kreiss_fig} (left) also shows large regions of pseudospectra corresponding to transient (off-limit cycle) behavior. \cref{DMD_kreiss_fig} (right) shows the Kreiss constant of the Koopman operator as $h$ varies.

The physical interpretation of these results aligns with those in \cref{sec:floquet} in that it is physically reasonable for the system to be stable while exhibiting large regions of pseudospectra. However, the Kreiss constants vary significantly between the methods and show a different type of dependency on $h$. We should not be alarmed by the differences in the Kreiss constants, as the norms they are based on are fundamentally different and incomparable. Furthermore, the Koopman operator provides a global linearization of a dynamical system, in contrast to the local linearization discussed in \cref{sec:floquet}.

\section{Conclusion}

Two fundamental questions in using pseudospectra for mathematical biology are: ``How should the matrix be obtained?" and ``How are clinically relevant perturbations modeled?". We addressed the first question in the context of the HPA axis by calculating pseudospectra through three methods: a time-dependent Jacobian, linearization around the limit cycle, and global linearization via Koopman operators. In addressing ``How should the matrix be obtained?", our initial model \cite{malek1, malek2} represents one of the many possible approaches. Future work could consider models incorporating CRH as a variable (as opposed to a parameter) or including glucocorticoid receptors. However, we want to highlight the role of spectra and pseudospectra in what we term ``model substantiation," assessing how well a model aligns with experimental evidence. Spectra and pseudospectra are useful for identifying the limits of or substantiating models based on their perturbation response. Other tests for model substantiation include assessing the stability of a limit cycle and the existence of oscillations under physically reasonable parameters.

Our first method, the time-dependent Jacobian, provided physically relevant results consistent with early rat experiments \cite{rats1, rats2}, where aggressive behavior was observed when rats were perturbed on the upward slope of corticosterone secretion. Additionally, fluctuations in the non-normality index indicated that the rhythms of ACTH and cortisol could produce a range of responses within a single period. Linearization around the limit cycle showed that, although the cycle was stable, the response to perturbations could be disproportionate, as evidenced by the pseudospectra. The DMD method similarly demonstrated stability across all eigenvalues and limit cycle-like structures in the eigenvectors yet highlighted sensitivity to perturbations in the pseudospectra. For computing pseudospectra in the latter two cases, we developed novel mathematical techniques for Banach spaces and for DMD in delay differential equations. Interestingly, increasing $h$ resulted in diminishing transient effects in local linearization while enhancing them in global linearization, as captured by the Kreiss constants.

Future work will focus on modeling clinically relevant perturbations. In \cite{rats2}, a male rat faced a stressor (a competitor rat) during both the rising and falling phases of cortisol. In \cite{rats1}, adrenalectomized rats received either a constant or a pulsatile dose of corticosterone; while the ACTH response was dampened with a constant dose, oscillatory corticosterone plasma levels and corresponding aggressive behavior were still observed. Our current models and pseudospectra methods do not account for extended perturbation durations or differences between injected and plasma corticosterone. To refine our understanding of aggressive behavior during unstable periods, as observed in \cref{default_h_fig}, we propose using structured pseudospectra -- specifying the $E$ in the pseudospectra definition -- to better replicate experimental outcomes. Additionally, linearization around the limit cycle and DMD could elucidate the effects of hormonal misalignment, exacerbated by external stressors like sleep deprivation and meal timing \cite{zavala2022misaligned}. Converting these stressors into structured ACTH perturbations will be a key focus.

Finally, data-driven approaches like DMD help bypass challenges in model substantiation but introduce difficulties in distinguishing perturbations from underlying rhythms or deriving perturbations from data. With established prior rhythms and modelable perturbations, we can investigate how structured pseudospectra on matrices align models with reality. This opens avenues for studying individual variability in perturbation responses and advancing personalized medicine. We advocate for strong collaboration between mathematicians and experimentalists to ensure data suits these sophisticated techniques.


\vspace{1cm}

\begingroup
\let\clearpage\relax
\bibliography{HPA_axis.bib}
\endgroup

\end{document}